\documentclass{article}
\usepackage{amsmath,amsfonts,amssymb,amscd}
\usepackage[T2A]{fontenc}
\usepackage{graphicx}

\newtheorem{theo}{Theorem}
\newtheorem{lem}{Lemma}[section]
\newtheorem{prop}{Proposition}[section]
\newtheorem{cor}{Corollary}[section]
\newtheorem{rem}{Remark}[section]

\makeatletter \@addtoreset{equation}{section} \makeatother

\newcommand{\mC}{\mathbb{C}}

\newcommand{\mR}{\mathbb{R}}

\newcommand{\mZ}{\mathbb{Z}}
\newcommand{\mN}{\mathbb{N}}

\newcommand{\ba}{{\bf a}}
\newcommand{\bh}{{\bf h}}
\newcommand{\bs}{{\bf s}}

\newcommand{\AAA}{{\cal A}}

\newcommand{\CC}{{\cal C}}

\newcommand{\EE}{{\cal E}}
\newcommand{\FF}{{\cal F}}
\newcommand{\HH}{{\cal H}}
\newcommand{\LL}{{\cal L}}
\newcommand{\MM}{{\cal M}}

\newcommand{\eps}{\varepsilon}
\newcommand{\ph}{\varphi}
\newcommand{\thet}{\vartheta}

\newcommand{\ssupp}{\operatorname{sing\,supp}}
\newcommand{\supp}{\operatorname{supp}}

\newcommand\qed{{\unskip\nobreak\hfil\penalty50
  \hskip2em\hbox{}\nobreak\hfil\mbox{\rule{1ex}{1ex} \qquad}
    \parfillskip=0pt \finalhyphendemerits=0\par\medskip}}

\begin{document}

\title
{Oscillator and thermostat}
\author{D.Treschev}
\maketitle

\begin{abstract}
We study the problem of a potential interaction of a
finite-dimensional Lagrangian system (an oscillator) with a linear
infinite-dimensional one (a thermostat). In spite of the energy
preservation and the Lagrangian (Hamiltonian) nature of the total
system, under some natural assumptions the final dynamics of the
finite-dimensional component turns out to be simple while the
thermostat produces an effective dissipation.
\end{abstract}

\section{Introduction}

The problem of interaction of a finite-dimensional Hamiltonian
system with an infinite-dimensional one always attracted attention
of physicists and specialists in dynamics. In comparison with the
finite-dimensional dynamics in spite of the energy preservation
in this situation in general some effective dissipation appears.
More precisely, a part of the energy can tend to be ``unobservable''
because of a radiation (existence of waves, going to infinity) and
because of transformation into heat (when some $L_2$-component of
a solution oscillates faster and faster, and weakly tends to zero,
keeping a noticeable part of the total energy). Note that by using
the Fourier transform one can see that these two mechanisms are dual.

The problems of this kind are complicated and usually the
infinite-dimensional system has to be chosen linear while the
finite-dimensional system can be linear or nonlinear. The
infinite-dimensional part can be a thermostat \cite{Bog}, a string
(or more generally, the Klein-Gordon system) \cite{Kom,KSK97},
or an electromagnetic field (the Maxwell equations) \cite{KS00}.

Below we regard the infinite-dimensional system as a continual collection of independent oscillators. We call it a thermostat. Note that thermostat is only one of possible physical interpretations of the infinite-dimensional system in our context, probably, not quite unquestionable.

In this paper we assume that total energy of the system is finite. Our method is based on the simple observation that the finite-dimensional system should oscillate in a way which does not produce a resonance in the thermostat. (Existence of such a resonance implies infiniteness if the energy.) This condition imposes some restrictions on the Fourier spectrum of the finite-dimensional component of the solution. In the examples (Section \ref{sec:exa}) these restrictions imply simple final dynamics of the finite-dimensional subsystem.

In Section \ref{sec:exa} we consider a one-dimensional oscillator. In the nonlinear case we present some natural sufficient conditions for the thermostat such that the oscillator tends to equilibrium positions as $t\to\pm\infty$. If the oscillator is linear, some kind of synchronization can happen and the final dynamics of the system can be a harmonic oscillation with some ``eigenfrequency''. In both cases a part of the total energy ``dissolves'' in the thermostat, i.e., it is concentrated in a component of the thermostat motion which weakly tends to zero.

Now we turn to more technical part.
Consider the Lagrangian system
\begin{equation}
\label{L0}
     \frac d{dt} \frac{\partial \LL_0}{\partial \dot x}
   - \frac{\partial \LL_0}{\partial x}
 = 0, \qquad
   \LL_0
 = \LL_0(x,\dot x), \quad
   x\in M.
\end{equation}
Here $M$ is a smooth $m$-dimensional manifold,
$x = (x_1,\ldots,x_m)$ are local coordinates, and the Lagrangian
$\LL_0$ equals the difference of the kinetic and potential energy
$$
     \LL_0(x,\dot x)
  = \frac12 \langle A(x) \dot x,\dot x\rangle - V_0(x).
$$
The square $m\times m$-matrix $A$ is assumed to be positive definite.
Hence, it determines a Riemannian metric on $M$. Then the product
$\langle A(x) \dot x,\dot x\rangle$ is the square of the Riemannian
length of the velocity vector $\dot x$. Such Lagrangian systems are
usually called {\it natural}.

Consider another system, called a thermostat. This is an
infinite-dimensional linear Lagrangian system with Lagrangian
\begin{equation}
\label{L_T}
    \LL_T(\xi,\dot\xi)
  = \int \frac {\rho(\nu)}{2}
         \big( \dot\xi^2(\nu,t) - \nu^2\xi^2(\nu,t)
         \big)\,d\nu .
\end{equation}
(For brevity we use the notation $\int$ for $\int_{-\infty}^{+\infty}$.)
One can regard it as a continual collection of independent harmonic oscillators, parameterized by the internal frequency of oscillations $\nu$. Physical meaning of $\rho(\nu)$ is the density of oscillators with the frequency $\nu$. Equations of motion of the isolated thermostat are
$$
    \rho(\nu) \big( \ddot\xi(\nu,t) + \nu^2 \xi(\nu,t)
              \big)
  = 0.
$$

Now suppose that these two systems interact and the interaction
potential is linear in $\xi$:
$$
    V_{\mbox{\rm int}}
  = - f(x)\, \int \kappa(\nu)\,\xi\, d\nu,
$$
where $f : M\to\mR$ is a smooth function. In other words we
consider the system with Lagrangian
$\LL_0 + \LL_T - V_{\mbox{\rm int}}$. The equations of motion are as
follows:
\begin{eqnarray}
\label{1st}
     \frac d{dt} \frac{\partial \LL_0}{\partial\dot x}
           - \frac{\partial \LL_0}{\partial x}
 &=&   \frac{\partial f}{\partial x}\,
       \int \kappa\xi\, d\nu , \\
\label{2nd}
     \rho ( \ddot\xi + \nu^2 \xi)
 &=&  \kappa f(x) .
\end{eqnarray}
We take initial values in the form
\begin{equation}
\label{ini}
        x|_{t=0} = x_0, \quad
      \xi|_{t=0} = \xi_0(\nu), \quad
   \dot x|_{t=0} = \dot x_0, \quad
  \dot\xi|_{t=0} = \dot\xi_0(\nu).
\end{equation}
The total energy is
\begin{equation}
\label{ene}
     E(x,\xi,\dot x,\dot\xi)
  =  \frac12 \langle A\dot x,\dot x\rangle
    + \frac12\int \rho\dot\xi^2\, d\nu
    + V_0 + \frac12 \int \rho\nu^2\xi^2\, d\nu
    + V_{\mbox{int}}.
\end{equation}

The paper is organized as follows. Section \ref{sec:mth} starts from formulation of our assumptions on the functions $\xi_0,\dot\xi_0,\rho,\kappa,A,V_0$, and $f$. A part of these assumptions is responsible for regularity of these functions. Another part ensures positive definiteness of the energy $E$. We need this positive definiteness to prove Theorem \ref{theo:existence} on the existence of a solution to system (\ref{1st}),(\ref{2nd}),(\ref{ini}).

We consider the space $\CC_b$ of bounded uniformly continuous functions on a line and the space $\hat\CC_b$ of distributions, obtained as
Fourier transforms of functions from $\CC_b$. The subspace $\CC_0\subset\CC_b$ consists of the functions, vanishing at infinity. For any distribution $\hat\psi\in\hat\CC_b$ we define its singular support $\ssupp\hat\psi\subset\mR$ as follows. We say that $\lambda\in\ssupp\hat\psi$ if there exists a regular (for example, smooth) function $\hat h$ with $\supp\hat h$ in an arbitrarily small neighborhood of $\lambda$ such that Fourier transform of
$\hat h\hat\psi$ does not lie in $\CC_0$.

Except Theorem \ref{theo:existence} Section \ref{sec:mth} contains formulation of Theorem \ref{theo:exotic}, an important technical tool.
Assertion (1) of Theorem \ref{theo:exotic} means that the Fourier transform $\hat\psi$ of the function $\psi(t) = f(x(t))\in\CC_b$ is regular on the ``essential spectrum'' $\sigma$ of the thermostat. Here $\sigma\subset\mR$ is the set of essential frequencies of oscillators in the thermostat.
The physical meaning of $\sigma$ is as follows.
Oscillators from the thermostat with frequencies outside of $\sigma$
do not influence the finite-dimensional part of the system and therefore, can be ignored. Formally, $\sigma = \{\nu : \hat a(\nu) \ne 0\}$,
$\; \hat a = \kappa / \sqrt\rho$. Regularity of $\hat\psi$ on $\sigma$ means that $\ssupp\hat\psi \cap \sigma = \emptyset$. Assertion (2) of Theorem \ref{theo:exotic} can be regarded as a relation between the singular part of $\hat\psi$ and the singular part of $\hat\ph$, the Fourier transform of $\int\kappa(\nu)\xi(\nu,t)\,d\nu$. The last integral is taken from the right-hand side of (\ref{1st}).

Section \ref{sec:exa} contains applications of our technical observations to the case of one-dimensional oscillator. First we show that the well-known Klein-Gordon equation generates a thermostat, i.e.,
it is equivalent to an infinite-dimensional Lagrangian system with Lagrangian (\ref{L_T}). We define thermostats used in this section as the ones, satisfying Hypothesis {\bf T}. Conditions, presented in it are motivated by the properties of the Klein-Gordon thermostat.

Then we specify the finite-dimensional part of the system, by considering a one-dimensional oscillator with
$\LL_0 = \frac12 \dot x^2 - V_0(x)$ with some smooth potential $V_0$. We take $f(x)=x$ and distinguish two cases: linear ($V_0=vx^2/2$) and nonlinear ($V_0$ is arbitrary).

In the linear case we show that, when $t\to\pm\infty$, depending on
parameters of the system the oscillator tends to the equilibrium
position (Theorem \ref{theo:empty}) or to some harmonic motion with
an eigenfrequency $\lambda_0$ (Theorem \ref{theo:linear}). The last
situation can be regarded as some kind of synchronization in the total
system.
In the nonlinear case we have to assume that the spectrum $\sigma$ is very large ($\sigma = \mR\setminus\{0\}$). Then we obtain that the oscillator tends to some equilibrium positions which can be different for
$t\to +\infty$ and $t\to-\infty$, and belong to the set of critical points of the effective potential $V_0(x) - \frac12 Kx^2$,
$K = \int \frac{\kappa^2}{\rho\nu^2}\, d\nu$.

In Section \ref{sec:existence} we obtain some a priori estimates which follow from the energy preservation and prove Theorem \ref{theo:existence}. In Section \ref{sec:psi} we show that preservation of the energy implies that the $L_2$-norm of the functions
$\hat a\hat\psi_t$ is uniformly bounded (Corollary \ref{cor:rhoxi}).
Here $\hat\psi_t(\tau)$ is the Fourier transform of
$\psi(\tau) \chi_{(0,t)}(\tau)$, where $\chi_{(0,t)}$ is the characteristic function of the interval $(0,t)$. This statement is crucial in the proof of assertion (1) of Theorem \ref{theo:exotic}
(Section \ref{sec:lemmas}). Assertion (2) of Theorem \ref{theo:exotic} is proven in Section \ref{sec:kappaxi}. Section \ref{sec:CCM} contains a self-contained theory of the spaces $\CC_b,\CC_0$ and others together with their Fourier transforms.

\section{Technical theorems}
\label{sec:mth}

Below the functions $\xi_0,\dot\xi_0$, $\rho,\kappa$, $A,V_0$, and $f$ satisfy hypotheses {\bf H1}--{\bf H5}.
\medskip

{\bf H1}. {\it The functions $A$, $V_0$, $f$ are smooth, $A$ is positive definite, $\rho$ is non-negative, and $\supp\kappa\subset\supp\rho$.}
\smallskip

Below we extend the function $\kappa / \sqrt\rho$ to $\mR\setminus\supp\rho$ as zero.

{\bf H2}. {\it The integrals
$\displaystyle{\int \frac{\kappa^2}{\rho}\,d\nu}$ and
$\displaystyle{K := \int \frac{\kappa^2}{\rho\nu^2}\,d\nu}$
are finite.}

{\bf H3}. {\it The effective potential
$$
  V(x) = V_0(x) - K f^2(x) / 2
$$
is bounded from below (without loss of generality $V\ge 0$). Moreover,
either $f$ is bounded or $V\to\infty$ as $|f|\to\infty$. The sets
$M_E = \{x\in M : V\le E\}$ are compact.}
\medskip

In particular, by (\ref{ene})
\begin{equation}
\label{leE}
     E
  =  \frac12 \langle A\dot x,\dot x\rangle
   + \frac12 \int \rho \dot\xi^2\, d\nu
   +  V
   +  \frac12 \int \rho\nu^2\Big( \xi - \frac{\kappa f}{\rho\nu^2}
                            \Big)^2 \, d\nu,
\end{equation}
where all four terms in the right-hand side are non-negative.

{\bf H4}. $\int \rho\, (\dot\xi_0^2 + \nu^2\xi_0^2)\, d\nu < \infty$.

\begin{rem}
\label{rem:E<infty}
Assumptions {\bf H2} and {\bf H4} imply that
$E(x_0,\xi_0,\dot x_0,\dot\xi_0) < \infty$. Indeed, by (\ref{ene}) it is sufficient to check that $\int \kappa\xi_0\,d\nu < \infty$. This inequality follows from the estimate
$\big( \int \kappa\xi_0\,d\nu \big)^2
  \le K \int\rho\nu^2\xi_0^2\,d\nu$.
\end{rem}

\begin{theo}
\label{theo:existence}
Suppose that conditions {\bf H1}--{\bf H4} hold. Then there exists a solution $x(t),\xi(\nu,t)$ of equations (\ref{1st}), (\ref{2nd}), (\ref{ini}). Moreover, $E<\infty$ is constant and the functions
\begin{equation}
\label{psiE}
    \psi(t)
  = f(x(t)), \quad
    \EE(t)
  = \int \rho(\nu) \Big( \nu^2\xi^2(\nu,t) + \dot\xi^2(\nu,t)
                   \Big)\, d\nu
\end{equation}
are bounded.
\end{theo}

Proof of Theorem \ref{theo:existence} is contained in Section \ref{sec:existence}.
\medskip

Let $\CC_b$ be the space of uniformly continuous bounded functions $\mR\to\mC$. It is a Banach space with the $L_\infty$-norm:
$\|\ph\|_\infty = \sup |\ph|$. Let $\MM$ be the space of densities of finite Radon $\mC$-measures on $\mR$. The space $\MM$ is conjugated to $\CC_b$ and consists of distributions such that for any $\mu\in\MM$
$\,\mu(\tau)d\tau$ is a finite Radon $\mC$-measure. Therefore $\MM$ is a Banach space with respect to the norm
$$
    \|\mu\|_{\MM}
  = \sup_{\|\ph\|_\infty = 1} \mu(\ph), \qquad
    \mu(\ph) = \int \ph(\tau)\mu(\tau) \, d\tau, \quad
    \ph\in\CC_b, \; \mu\in\MM.
$$
The space $\MM$ obviously contains $L_1 = L_1(\mR)$.

Let $\FF$ be the Fourier transform:
$$
     \FF(\psi)
  =  \hat\psi(\lambda)
  =  \int e^{-i\lambda \tau} \psi(\tau)\, d\tau , \quad
     \FF^{-1}(\hat\psi)
  =  \psi(\tau)
  =  \frac 1{2\pi} \int e^{i\lambda \tau} \hat\psi(\lambda)\, d\lambda .
$$
Recall the standard formulas for the convolution
$$
      2\pi\,\FF(\ph\psi)
   = \FF \ph * \FF \psi, \quad
      \FF^{-1}(\hat\ph\hat\psi)
   = \FF^{-1} \hat\ph * \FF^{-1} \hat\psi.
$$

Below we usually denote functions (or distributions) on the time axis by $h,\ph,\psi,\ldots$ and their Fourier transforms by
$\hat h,\hat\ph,\hat\psi,\ldots$ Analogously, a hat over a notation of a function space denotes its Fourier dual space.

Let $\CC_0^+\subset\CC_b$ denote the space of functions which tend to zero at $+\infty$ and $\CC_0^-\subset\CC_b$ the space of functions which tend to zero at $-\infty$. We have the spaces of distributions
$\hat\CC_b = \FF(\CC_b)$ and $\hat\CC_0^\pm = \FF(\CC_0^\pm)$.
We put $\CC_0 = \CC_0^+\cap\CC_0^-$ and $\hat\CC_0 = \FF(\hat\CC_0)$.

We define $\hat\MM = \FF(\MM)$. It is easy to show that $\hat\MM\subset\CC_b$, \cite{Boga}. By the Riemann-Lebesgue theorem $\FF(L_1) \subset \CC_0$.

We put
\begin{equation}
\label{tildehat}
    \hat a(\nu)
  = \frac{\kappa(\nu)}{\sqrt{\rho(\nu)}},\quad
    \hat w(\nu)
  = \frac{2\pi\hat a^2(\nu)}{i\nu} .
\end{equation}

{\bf H5}. {\it The functions $\hat a$ and $\hat a' = d\hat a / d\nu$ belong to the space $\hat\MM$ and
$\hat a(0) = \hat w(0) = 0$.}
\medskip

Note that by Lemma \ref{lem:a-w} $w = \FF^{-1}(\hat w) \in L_1$.
\medskip

We also define the ``essential spectrum'' of the thermostat
$$
    \sigma
  = \{s : \hat a(s) \ne 0\}
  = \{s : \hat w(s) \ne 0\}.
$$
The motion of the thermostat outside the essential spectrum $\xi|_{\mR\setminus\sigma}$ is immaterial for the finite-dimensional part of the system because the coupling term in (\ref{1st}),
$\frac{\partial f}{\partial x} \int \kappa\xi\, d\nu$, is independent of
$\xi|_{\mR\setminus\sigma}$.

By {\bf H5} $\hat a$ is continuous. Therefore the set $\sigma$ is open and its closure $\overline\sigma = \supp\hat a = \supp\hat w$.
\medskip

Let $\hat\ph\in\hat\CC_b$, $\lambda\in\mR$. We say that $\lambda\not\in\ssupp^+(\hat\ph)$ iff there exists an interval
$I\ni\lambda$ such that for any $\hat\mu\in\hat\MM$,
$\supp\hat\mu \subset I$ we have: $\hat\mu\hat\ph\in \hat\CC_0^+$.
Analogously we define $\ssupp^-(\hat\ph)$ and
$\ssupp(\hat\ph) = \ssupp^+(\hat\ph) \cup \ssupp^-(\hat\ph)$.

Hence, $\ssupp^\pm(\hat\ph)$ and $\ssupp(\hat\ph)$ are closed subsets of $\mR$.
\medskip

We put
\begin{equation}
\label{wdia}
   w_\diamond^\pm(\tau)
 = \left\{ \begin{array}{cl}
        \pm (w(\tau) - w(-\tau)) / 2 \quad & \mbox{if } \pm\tau > 0, \\
                   0                 \quad & \mbox{if } \pm\tau < 0.
            \end{array}
     \right.
\end{equation}
A standard computation shows that\footnote{
  see some details in Section \ref{sec:kappaxi}}
\begin{equation}
\label{wdiamond}
    \hat w_\diamond^\pm(\nu)
  = \pm \frac14 \big( \hat w(\nu) - \hat w(-\nu) \big)
   + \frac 1{2\pi i}\,\, \mbox{v.p.}\!\!
      \int \frac{\lambda \, \hat w(\lambda)\, d\lambda}
                {\nu^2 - \lambda^2} .
\end{equation}

\begin{prop}
\label{prop:diamond}
Suppose that $\hat a$ and $\hat w$ satisfy {\bf H5}. Then

(a) $w_\diamond^\pm := \FF^{-1}(\hat w_\diamond^\pm) \in L_1$,

(b) $\mbox{\rm Im}\,\hat w_\diamond^\pm$ are odd and
$\mbox{\rm Re}\,\hat w_\diamond^\pm$ are even,

(c) $\hat w_\diamond^\pm(0)
 = -\frac{1}{2\pi i}\int \frac{\hat w(\lambda)}{\lambda}\, d\lambda
 = K$, the constant, defined in {\bf H2},

(d) $\hat w_\diamond^\pm$ are $C^1$-smooth.
\end{prop}

{\it Proof}. By {\bf H5} and Lemma \ref{lem:a-w} $w\in L_1$. Therefore by (\ref{wdia}) $w_\diamond^\pm$ also lie in $L_1$.

Assertions (b) and (c) are obvious while (d) follows from {\bf H5}. \qed

\begin{theo}
\label{theo:exotic}
Suppose that conditions {\bf H1}--{\bf H5} hold.
Then the functions $\psi(t) = f(x(t))$ and
$\ph(t) = \int\kappa(\nu)\,\xi(\nu,t)\,d\nu$ are bounded and their Fourier transforms $\hat\psi = \FF(\psi)$, $\hat\ph = \FF(\ph)$ are such that

(1) $\hat a\hat\psi\in L_2$, $\hat w\hat\psi\in L_1$, and
    $\ssupp \hat\psi \cap \sigma = \emptyset$,

(2) $\hat\ph = \hat w_\diamond^+ \hat\psi + \hat\thet^+
             = \hat w_\diamond^- \hat\psi + \hat\thet^-$ for some
   $\hat\thet^\pm\in\hat\CC_0^\pm$.
\end{theo}

We prove Theorem \ref{theo:exotic} in Sections \ref{sec:kappaxi}--\ref{sec:lemmas}.

\begin{rem}
By (\ref{1st})
\begin{equation}
\label{(2)}
     \displaystyle{
      \frac d{dt} \frac{\partial \LL_0}{\partial \dot x}
     - \frac{\partial \LL_0}{\partial x}
 =  \frac{\partial f}{\partial x} \, \ph. }
\end{equation}
\end{rem}

{\bf Some explanations}.
Physical meaning of assertion (1) from Theorem \ref{theo:exotic} is a regularity of $\hat\psi$ on $\sigma$. This follows essentially from the fact that if $\ssupp\hat\psi \cap\sigma\ne\emptyset$,
a resonance appears in the equation
$$
    \rho(\nu) \big( \ddot\xi + \nu^2 \xi \big)
  = \psi(t)
$$
Existence of such a resonance would contradict to the assumption that
the energy is bounded.
\smallskip

Assertion (2) + Remark 1 give a possibility to write down an approximate
($\bmod\,\CC_0^\pm$) equation of motion for the oscillator. Studying of a final dynamics of the oscillator is based on the analysis of this equation.

\section{Examples}
\label{sec:exa}

\subsection{Klein-Gordon thermostat}

The 1-dimensional Klein-Gordon equation has the form
$u_{tt} = u_{qq} - m_0^2 u$, $m_0\ge 0$. If the constant $m_0$ vanish, the Klein-Gordon equation coincides with the 1-dimensional wave equation (a string). Passing to the Fourier transform
$\xi(s,t) = \int e^{-iqs} u(q,t)\, dq$, we obtain the equation
$\ddot\xi = -(s^2 + m_0^2)\xi$. The corresponding Lagrangian is
$\LL_T
 = \frac12 \int \rho_0(s)\big(\dot\xi^2 - (s^2 + m_0^2)\xi^2\big)
   \, ds$ with an arbitrary function $\rho_0$.
Putting $\nu^2 = m_0^2 + s^2$, $\nu s \ge 0$, we obtain:
$$
    \LL_T
  = \frac12 \int \rho(\nu)(\dot\xi^2 - \nu^2\xi^2)\, d\nu, \qquad
    \rho(\nu)
  = \left\{  \begin{array}{cl}
                0       &   \mbox{ if } |\nu| \le m_0 , \\
             \frac{\rho_0\big( \nu \sqrt{1 - m_0^2\nu^{-2}}\big)}
                  {\sqrt{1 - m_0^2\nu^{-2}}}
                        &   \mbox{ if } |\nu| \ge m_0 .
             \end{array}
    \right.
$$
Below this system is called a Klein-Gordon thermostat.
We have the following obvious

\begin{prop}
\label{prop:KG}
Suppose that in a Klein-Gordon thermostat $\rho_0(s)$ is $C^1$-smooth and positive for $s\ne 0$, $\lim_{s\to 0} s / \rho_0(s) < \infty$, $\kappa$ is $C^1$-smooth, positive, and $\supp\kappa = \mR\setminus (-m_0,m_0)$. Then

(a) $\rho$ is continuous and positive for $|\nu| > m_0$;

(b) the functions
$\hat a$ and $\hat w$ (see (\ref{tildehat})) vanish for $|\nu|\le m_0$, positive for $|\nu| > m_0$, and $C^1$-smooth.

(c) the set $\sigma$ is
\begin{equation}
\label{sigmaKG}
  \sigma = \{\nu\in\mR : |\nu|>m_0 \}.
\end{equation}

(d) the functions $\hat w_\diamond^\pm$ (see (\ref{wdiamond})),
restricted to the interval $[-m_0,m_0]$, equal
$$
    \hat w_\diamond^\pm(\nu)
  = \int\frac{\hat a^2(\lambda)}{\lambda^2 - \nu^2} \, d\lambda, \qquad
    |\nu| \le m_0,
$$
i.e., they are real, even, and increasing for $0\le\nu\le m_0$.
\end{prop}

Below we consider a thermostat, satisfying the following hypothesis, motivated by Proposition \ref{prop:KG}.
\smallskip

{\bf T}. (a) $\hat a(\nu)$ vanishes for $|\nu|\le\nu_0$ and positive
  for $|\nu| > \nu_0$,

$\quad\,$ (b) $\sigma = \{\nu\in\mR : |\nu| > \nu_0 \}$,

$\quad\,$ (c) $\hat w^\pm_\diamond|_{[-\nu_0,\nu_0]}$ are equal to each other, real, even and increasing on $[0,\nu_0]$.

\subsection{One-dimensional oscillator}

Consider the system (\ref{1st}),(\ref{2nd}),(\ref{ini}) where
$$
  x\in\mR, \quad
  A = 1, \quad
  f = x.
$$
Then $\psi(t) = x(t)$ and equation (\ref{(2)}) takes the form
\begin{equation}
\label{[L]=}
  \ddot\psi + V'_0(\psi) = \ph, \qquad
  \hat\ph = \hat w_\diamond^+ \hat\psi + \hat\thet^+
          = \hat w_\diamond^- \hat\psi + \hat\thet^-, \qquad
  \hat\thet^\pm\in\hat\CC_0^\pm.
\end{equation}

In the next two sections we consider the dynamics in the system oscillator + thermostat in two cases.

(A) $V_0 = vx^2 / 2$, $v > 0$.

(B) $\nu_0 = 0$.

\subsection{Case A (linear)}

Suppose that $V_0 = vx^2 / 2$, $v > 0$. By {\bf H3} we have to assume that $v - K > 0$. We define
$$
    \hat\Phi^\pm(\nu)
  = - \nu^2 + v - \hat w_\diamond^\pm(\nu).
$$
Note that by {\bf T} $\hat\Phi^\pm$, restricted to
$\mR\setminus\sigma = [-\nu_0,\nu_0]$, are real, even, and decreasing for
$0\le\nu\le\nu_0$. We have: $\Phi^\pm(0) = v - K > 0$.
The Fourier transform of (\ref{[L]=}) yields
\begin{equation}
\label{linear}
    \hat\Phi^\pm \hat\psi
  = \hat\thet^\pm, \qquad
    \hat\thet^\pm \in \hat\CC_0^\pm .
\end{equation}

\begin{rem}
\label{rem:ss}
By Theorem \ref{theo:exotic} and Lemma \ref{cor:1/Phi} equation (\ref{linear}) implies that
$\ssupp^\pm(\hat\psi) \subset \{\hat\Phi^\pm = 0\} \setminus \sigma$.
\end{rem}

\begin{theo}
\label{theo:empty}
Suppose that a linear one-dimensional oscillator interacts with a thermostat, satisfying ${\bf T}$, conditions {\bf H1}--{\bf H5} hold, and
$\hat w^\pm_\diamond(\nu_0) < v$. Then $\psi\in\CC_0$.
\end{theo}

{\it Proof}. We have: $\Phi^\pm(0) > 0$ and $\Phi^\pm(\nu_0) > 0$.
Hence by monotonicity and evenness $\Phi^\pm \ne 0$ on
$[-\nu_0,\nu_0] = \mR\setminus\sigma$. By Remark \ref{rem:ss}
$\ssupp^\pm(\hat\psi) = \emptyset$. Now $\psi\in\CC_0$ by Corollary \ref{cor:ssupp=0}. \qed

Theorem \ref{theo:empty} means that, if $\hat w^\pm_\diamond(\nu_0) < v$ on $\mR\setminus\sigma$, a linear one-dimensional oscillator transmits all its energy to a thermostat when $t\to\pm\infty$.
\smallskip

\begin{theo}
\label{theo:linear}
Suppose that a linear one-dimensional oscillator interacts with a thermostat, satisfying ${\bf T}$, conditions {\bf H1}--{\bf H5} hold, and
$\hat w^\pm_\diamond(\nu_0) > v$. Then for some $\alpha\in \mC$
\begin{equation}
\label{alphaalpha}
    \psi(t)
  = \alpha e^{i\lambda_0 t} + \overline\alpha e^{-i\lambda_0 t}
   + \thet(t), \qquad
    \thet \in \CC_0.
\end{equation}
\end{theo}

In other words, equation (\ref{alphaalpha}) implies that the final motion of the oscillator is reduced to harmonic oscillations with the frequency $\lambda_0$.

{\it Proof}. Putting $\eta = \sqrt\rho \xi$, we present equations (\ref{1st}),(\ref{2nd}),(\ref{ini}) in the form
\begin{eqnarray}
\label{xeta}
 &&\!\!\!\!
   \ddot x  = - vx + \int \hat a\eta \, d\nu,\quad
  \ddot\eta = -\nu^2\eta + \hat a x, \\
\label{xeta_ini}
 &&\!\!\!\!
     x|_{t=0} = x_0,\quad
  \eta|_{t=0} = \eta_0(\nu), \quad
   \dot x|_{t=0} = \dot x_0,\quad
 \dot\eta|_{t=0} = \dot\eta_0(\nu)  .
\end{eqnarray}

Consider the Hilbert space $\HH$ with elements $\zeta = (x,\eta)$ and scalar product
$   (\zeta_1,\zeta_2)
  = x_1 \overline x_2 + \int \eta_1 \overline\eta_2 \, d\nu.
$
Then system (\ref{xeta}) can be written as
\begin{equation}
\label{z'=Az}
  \ddot\zeta = - \AAA\zeta, \qquad
  \zeta|_{t=0} = \zeta_0, \quad
  \dot\zeta|_{t=0} = \dot\zeta_0,
\end{equation}
where $\AAA$ is the self-adjoint unbounded operator
$$
             \zeta
    =        (x,\eta)
 \mapsto \AAA \zeta
    =    \Big( vx - \int \hat a\eta \, d\nu,\,
                - \hat a x + \nu^2\eta \Big).
$$
Inequality $v - K > 0$ implies that $\AAA$ is positive definite. Indeed:
$$
     (\zeta,\AAA\zeta)
  =  \int \Big( \nu\eta - \frac{\hat a}{\nu} x \Big)^2 \, d\nu
    + ( v - K) x^2
  > 0  \quad
     \mbox{for any $\zeta\ne 0$} .
$$
The solution of (\ref{z'=Az}) is
\begin{equation}
\label{sol}
    \zeta(t)
  = \frac12 e^{i\sqrt\AAA\,t} \zeta_+
   + \frac12 e^{-i\sqrt\AAA\,t} \zeta_- ,\qquad
    \zeta_\pm
  =  \zeta_0 \mp i\AAA^{-1/2} \dot\zeta_0 .
\end{equation}
By the spectral theorem \cite{Riesz-SN} we have:
$\AAA = \int_0^\infty s\, dE_s$, where $\{E_s\}$ is the corresponding spectral family. Since
$e^{\pm i\sqrt\AAA\, t} = \int_0^\infty e^{\pm i\sqrt s \, t}\, dE_s$, we present (\ref{sol}) in the form
\begin{equation}
\label{zeta=}
  \zeta(t)
  = \frac12 \int_0^\infty e^{i\sqrt s\,t} \, d E_s \zeta_+
   + \frac12 \int_0^\infty e^{-i\sqrt s\,t} \, d E_s \zeta_- .
\end{equation}

Direct computation shows that $\lambda_0^2 \in [0,\nu_0^2)$ is an eigenvalue of $\AAA$: $\; \AAA\zeta_{\lambda_0} = \lambda_0^2\zeta_{\lambda_0}$ for some $\zeta_{\lambda_0} \ne 0$ iff
$$
  \hat\Phi^\pm(\lambda_0) = 0, \quad
    \zeta_{\lambda_0}
  = \Big( 1,\frac{\hat a(\nu)}{\nu^2 - \lambda_0^2} \Big).
$$
(Recall that by Hypothesis {\bf T}(c) on $[-\nu_0,\nu_0]$ we have:
$\hat\Phi^+ = \hat\Phi^-$). By {\bf T}(a) continuous spectrum of $\AAA$ is $[\nu_0^2,\infty)$. Hence the eigenvalue $\lambda_0^2$ is an isolated point of $\mbox{spec}\,(\AAA)$.

We put ${\bf 1} = (1,0)\in\HH$. Consider the change of the variable $\lambda = \sqrt s$ in the first integral (\ref{zeta=}) and
$\lambda = -\sqrt s$ in the second one. Then the $x$-component of solution (\ref{zeta=}) is
$$
    x(t)
  = \int e^{i\lambda t} \hat\omega(\lambda)\, d\lambda,
$$
where $\hat\omega$ satisfies the equation
$$
    \hat\omega(\lambda)\, d\lambda
  = \left\{
     \begin{array}{cc}
       \lambda^{-1}\,d({\bf 1}, E_{\lambda^2} \zeta_+),
           & \quad \mbox{if } \lambda > 0, \\
       \lambda^{-1}\,d({\bf 1}, E_{\lambda^2} \zeta_-)
           & \quad \mbox{if } \lambda < 0.
     \end{array}
    \right.
$$
The distribution $\hat\omega$ is well-defined because $0\not\in\mbox{spec}\,\AAA$. Moreover, $\hat w\in\hat\CC_b$ because by
Theorem \ref{theo:existence} $\psi(t) = x(t)$ is bounded.

Since $\lambda_0^2$ is an isolated point of $\mbox{spec}(\AAA)$, in a small neighborhood of the point $\lambda = \lambda_0$ the distribution $\hat\omega(\lambda)$ equals $\alpha \delta(\lambda - \lambda_0)$, $\alpha\in\mC$. Analogously
$\hat\omega(\lambda) = \beta \delta(\lambda + \lambda_0)$ near $-\lambda_0$. Since $x$ is real, $\beta = \overline\alpha$. Therefore
$$
     x(t)
  =  \psi(t)
  =  \alpha e^{i\lambda_0 t} + \overline\alpha e^{-i\lambda_0 t}
    + \int e^{i\lambda t} \hat\omega_0(\lambda)\, d\lambda,
$$
where
$\hat\omega_0(\lambda)
  =  \hat\omega(\lambda) - \alpha_+ \delta(\lambda - \lambda_0)
    - \alpha_- \delta(\lambda + \lambda_0)$
vanishes in a neighborhood of the points $\pm\lambda_0$.

Now recall that by Remark \ref{rem:ss} for any $\lambda\in\ssupp^\pm\hat\psi$ we have: $\lambda = \lambda_0$ or $\lambda = -\lambda_0$. This implies that $\ssupp^\pm(\hat\omega_0) = \emptyset$. Applying Corollary \ref{cor:ssupp=0}, we finish the proof. \qed

\subsection{Case B (nonlinear)}

Now consider the case $\nu_0 = 0$ i.e., $\sigma = \mR\setminus\{0\}$.
\begin{theo}
\label{theo:nonlin}
Suppose that a one-dimensional (in general, nonlinear) oscillator interacts with a thermostat, satisfying {\bf T}, conditions {\bf H1}--{\bf H5} hold, and $\nu_0 = 0$. Then there exist limits
$\lim_{t\to\pm\infty} \psi(t) = x_\pm$, where $x_\pm$ are critical points of $V_0(x) - \frac12 Kx^2$.
\end{theo}

{\it Proof}. Since $\nu_0 = 0$, by assertion (1) of Theorem \ref{theo:exotic} $\ssupp^\pm(\hat\psi)$ contains at most one point: $\lambda = 0$.

By Lemma \ref{cor:onepoint} for any $\star\in\{+,-\}$ we have:
$\psi = \psi_0^\star + \psi_1^\star$, where $\supp\hat\psi_1^\star$ lies in a small neighborhood of 0,
$$
  \hat\psi_0^\star\in\hat\CC_0^\star, \quad
  \hat\psi_1^\star\in\hat\CC_b, \quad
  \psi_1^\star \in C^\infty, \quad
  \dot\psi_1^\star
 = \frac{d}{d\tau}\psi_1^\star \in\CC_0^\star.
$$
Since $\supp\hat\psi_1^\star$ is compact, we also have:
$\ddot\psi_1^\star
 = \FF^{-1}(-\lambda^2\hat\psi_1^\star)\in\CC_0^\star$.
By Proposition \ref{prop:diamond}, assertion (d), $\hat w_\diamond^\star$ is $C^1$-smooth at zero. Therefore
$\hat w_\diamond^\star\hat\psi_1^\star
 = \hat w_\diamond^\star(0) \hat\psi_1^\star \bmod\hat C_0^\star$. This
implies
$$
    \hat\ph
  = \hat w_\diamond^\star\hat\psi + \hat\thet^\star
  = \hat w_\diamond^\star(0) \, \hat\psi_1^\star + \hat\thet^\star_1, \qquad
    \hat\thet^\star_1
 \in \hat\CC_0^\star.
$$
Recall that $\hat w_\diamond^+(0) = \hat w_\diamond^-(0) = K$.
Then (\ref{[L]=}) takes the form
$$
  V'_0(\psi_1^\star) = K\psi_1^\star \bmod\CC_0^\star.
$$
Therefore for large $|t|$, $\mbox{sgn}\,t = \star$, the function $\psi_1^\star(t)$ can take values only in an arbitrarily small neighborhood of the set $\mbox{Cr} := \{x : V'_0(x) - Kx = 0\}$. This implies the existence of the limits
$$
   \lim_{t\to-\infty}\psi_1^-(t)
 = \lim_{t\to-\infty}\psi(t)
 = x_- \quad
   \mbox{and}\quad
   \lim_{t\to+\infty}\psi_1^+(t)
 = \lim_{t\to+\infty}\psi(t)
 = x_+ ,
$$
where $x_\pm \in \mbox{Cr}$. \qed

\section{Proof of Theorem \ref{theo:existence}}
\label{sec:existence}

The plan of the proof is standard: a combination of a local existence theorem with the energy estimate (\ref{leE}). Recall that by Remark \ref{rem:E<infty} $E(x_0,\xi_0,\dot x_0,\dot\xi_0) < \infty$.

\subsection{A priori estimates}

\begin{prop}
\label{cor:bounded}
If $E$ is finite, $\psi(t) := f(x(t))$ is bounded:
$|\psi| \le c_1 = c_1(E)$.
\end{prop}

Indeed, by {\bf H3} if $f$ is unbounded, $V$ is also unbounded. Then by (\ref{leE}) $E=+\infty$. \qed

\begin{prop}
\label{cor:kappaxi}
If $E$ is finite,
$\big| \int \kappa\xi\, d\nu \big| \le c_2 = c_2(E)$.
\end{prop}

Indeed, by the Cauchy-Bunyakovsky inequality
$$
     \Big( \int \kappa\xi\, d\nu \Big)^2
   = \Big( \int \frac{\kappa}{\sqrt\rho\,\nu} \,
                \sqrt\rho\,\nu\xi\, d\nu \Big)^2
 \le K \,\int \rho\nu^2 \xi^2 \, d\nu.
$$
By {\bf H2} $K<\infty$ while by (\ref{leE})
$$
     \int \rho\nu^2 \xi^2 \, d\nu
 \le 2 \Big( E + f(x) \int \kappa\xi\, d\nu \Big).
$$
Therefore
$ \big( \int \kappa\xi\, d\nu \big)^2
  < 2K \big( E + \psi\int \kappa\xi\, d\nu
       \big)$
which is equivalent to
$$
     \Big( \int \kappa\xi\, d\nu - K\psi
     \Big)^2
  <  2K E + K^2 \psi^2 .
$$
By Proposition \ref{cor:bounded} $\psi$ is bounded. Hence
$\int \kappa\xi\, d\nu$ is also bounded. \qed

\begin{prop}
\label{cor:ene-ineq}
If $E$ is finite, then $\EE$, defined in (\ref{psiE}), satisfies
$\EE \le c_3 = c_3(E)$.
\end{prop}

Indeed, by (\ref{leE})
$\EE - f\int\kappa\xi\,d\nu + \frac12 Kf^2 \le E$. It remains to use Propositions \ref{cor:bounded} and \ref{cor:kappaxi}. \qed

\subsection{Local existence theorem}

By (\ref{leE}) we have: $V \le E$. By {\bf H3} this implies that $x$ lies in a compact set $M_E\subset M$. It is possible to cover $M_E$ by a finite number of coordinate charts. Hence while proving a local existence theorem for solution with energy $E$ for the system (\ref{1st})--(\ref{2nd}), we can assume that $x\in B\subset\mR^m$, where
$$
  B = \{x\in\mR^m : \langle x,x\rangle \le 1\}
$$
is the unit ball with the center at zero. Moreover, we can assume that the parts of the charts, satisfying the inequality
$\langle x,x\rangle < 1/2$, also cover $M_E$, and
$$
  A^{-1}(x), f(x),\,
  \frac{\partial A}{\partial x_j},\,
  \frac{\partial f}{\partial x_j},\,
  \frac{\partial V_0}{\partial x_j}, \quad
  j = 1,\ldots,m \quad
  \mbox{are uniformly bounded and uniformly Lipschitz}.
$$

Putting
$$
    \xi_\pm
  = e^{\mp i\nu t} (\dot\xi \pm i\nu\xi), \quad
    p
  = A\dot x ,
$$
we represent system (\ref{1st})--(\ref{2nd}) in the form
\begin{equation}
\label{hamilt}
\begin{array}{l}
\displaystyle{
      \dot\xi_\pm \,
  =\, e^{\mp i\nu t} \, \frac{\kappa f(x)}{\rho}  , \quad
      \dot x \,
  =\, A^{-1} p,  }  \\
\displaystyle{
      \dot p
  =   \frac12 \big\langle G p,p
              \big\rangle
     - V_0'
     + f'\,\int \frac{\kappa}{2i\nu} \,
                \big( e^{i\nu t}\xi_+ - e^{-i\nu t}\xi_- \big)
                 \, d\nu , } \\
\end{array}
\end{equation}
where $G(x) = A^{-1}(x) A'(x) A^{-1}(x)$, $\xi_+ = \overline\xi_-$ and primes stand for $\partial / \partial x$.

This system is defined in an open set $U$ in the Hilbert space $\HH$ with elements $z = (\xi_+,\xi_-,x,p)$, where $x,p\in\mR^m$ and $\xi_+ = \overline\xi_-$. The scalar product is
$$
    (z,z)
  =      \langle A^{-1} p,p\rangle
   +  \langle x,x\rangle
   +  \frac12 \int \rho\, (|\xi_+|^2 + |\xi_-|^2)\, d\nu.
$$
We will also use the corresponding norm $\|z\|^2 :=(z,z)$.

Note that
\begin{equation}
\label{zz-E}
     \frac12 (z,z) - E
  =  \langle x,x\rangle - V + f \int\kappa\xi\,d\nu
     - Kf^2 / 2.
\end{equation}
The terms $-V$ and $- Kf^2 / 2$ are non-positive, and by Propositions
\ref{cor:bounded}--\ref{cor:kappaxi} $\; f\!\int\kappa\xi\,d\nu$ is bounded. Therefore for $x\in B$ (\ref{zz-E}) implies that $(z,z)$ is bounded. In other words, we can assume that $U\subset\HH$ is a bounded set.

The system (\ref{hamilt}) has the form $\dot z = \Lambda(z)$.

\begin{prop}
\label{prop:Lipsch}
The map $\Lambda : U\to\HH$ is Lipschitz.
\end{prop}

{\it Proof}. Take any two points
$z_1 = (\xi_{1+},\xi_{1-},x_1,p_1)$ and
$z_2 = (\xi_{2+},\xi_{2-},x_2,p_2)$ from $U$.

We want to prove that
$\|\Lambda(z_1) - \Lambda(z_2)\| \le C \|z_1 - z_2\|$. We have:
$\|\Lambda(z_1) - \Lambda(z_2)\| \le \sum_{\alpha=1}^5 \Delta_\alpha$, where
\begin{eqnarray*}
     \Delta_1^2
 &=& \int \rho \Big( \frac{\kappa f(x_1)}{\rho}
                    - \frac{\kappa f(x_2)}{\rho} \Big)^2\, d\nu , \\
     \Delta_2^2
 &=& \big\langle A^{-1} (p_1 - p_2), A^{-1} (p_1 - p_2) \big\rangle, \\
     \Delta_3^2
 &=& 3\big\langle A^{-1}(Q_1 - Q_2), Q_1 - Q_2 \big\rangle, \qquad
     Q_j
  =  \frac12 \big\langle G(x_j) p_j , p_j
             \big\rangle , \\
     \Delta_4^2
 &=& 3\big\langle A^{-1}(V_0'(x_1) - V_0'(x_2)),
                   V_0'(x_1) - V_0'(x_2) \big\rangle , \\
     \Delta_5^2
 &=& 3\big\langle A^{-1}(P_1 - P_2), P_1 - P_2 \big\rangle , \qquad
     P_j
  =   f'(x_j)\,\int \frac{\kappa}{2i\nu} \,
                    \big( e^{i\nu t}\xi_{j+} - e^{-i\nu t}\xi_{j-}
                    \big) \, d\nu .
\end{eqnarray*}

Since $f$ is Lipschitz,
$|f(x_1) - f(x_2)| \le \mbox{const}\, |x_1 - x_2|$. Therefore by
{\bf H2} for some $c,c_1>0$
$$
      \Delta_1^2
  \le c \int \frac{\kappa^2}{\rho}\, d\nu\, \|z_1 - z_2\|^2
  \le c_1\, \|z_1 - z_2\|^2  .
$$
Since $A^{-1}$ is uniformly bounded then
$\Delta_2^2 \le c_2 \|z_1 - z_2\|^2$.

Analogously we have:
$$
  \Delta_3^2 \le c_3 \|z_1 - z_2\|^2, \quad
  \Delta_4^2 \le c_4 \|z_1 - z_2\|^2.
$$
Here are some more details concerning $\Delta_5$:
$$
       \Delta_5^2
  \le  c \|z_1 - z_2\|^2\,
       \Big( \int \frac{|\kappa|}{\sqrt\rho\,|\nu|} \,
                   \sqrt\rho\, |\xi_+|\, d\nu
       \Big)^2
      + c\sup_{x\in B} |f'(x)|^2
        \Big(\int \frac{|\kappa|}{\sqrt\rho\,|\nu|} \,
            \sqrt\rho\, |\xi_{1+} - \xi_{2+}|\, d\nu
        \Big)^2
  \le  c_5 \, \|z_1 - z_2\|^2.
$$
\qed

By proposition \ref{prop:Lipsch} we can use the theorem on local existence of a solution for a system of ordinary differential equations \cite{Schwa}. Hence for any initial condition with $x$-component in the ball
$\{ \langle x,x\rangle < 1/2 \}$ we have a solution on the time interval $[-t_0,t_0]$ such that its $x$-component does not leave the ball
$\{ \langle x,x\rangle < 1 \}$. Here $t_0 > 0$ can be chosen the same for all on the charts. Since the energy preserves, we can continue the solution to all the time axis.

\section{Spectrum of $\psi(t)$}
\label{sec:psi}

For any measurable set $M\subset\mR$ let $\chi_M$ be its characteristic
function. If $M$ is an interval $(b_1,b_2)\subset\mR$, we take into account its orientation:
$$
    \chi_{(b_1,b_2)}(t)
  = \left\{ \begin{array}{cl}
             1 & \quad \mbox{if $\; b_1<t<b_2$}, \\
            -1 & \quad \mbox{if $\; b_2<t<b_1$}, \\
             0 & \quad \mbox{if $\; (b_1-t)(t-b_2) \le 0$}.
            \end{array}
    \right.
$$

Let $(x(t),\xi(\nu,t))$ be a finite energy solution of
(\ref{1st}),(\ref{2nd}),(\ref{ini}). For any real $t$ we put
\begin{equation}
\label{box}
    \psi_t
  = \chi_{[0,t]}\,\psi , \quad
    \psi(\tau)
  = f(x(\tau)) ,  \quad
     \hat\psi_t
  =  \FF(\psi_t) .
\end{equation}
Then $\hat\psi_t$ is an entire analytic function.
Since $\psi_t$ is real, we have:
\begin{equation}
\label{real}
  \hat\psi_t(\lambda)  =  \overline{\hat\psi}_t(-\lambda).
\end{equation}
Now we consider equation (\ref{2nd}) as if the function $\psi(t)$ were known. By (\ref{ini}) and (\ref{box}) we get:
$\xi(\nu,t) = \xi_1(\nu,t) + \xi_2(\nu,t)$, where
\begin{eqnarray}
\nonumber
      \xi_2(\nu,t)
 &:=& \xi_0(\nu) \cos \nu t
         + \frac{1}{\nu}\,\dot\xi_0(\nu) \sin \nu t
   =  \mbox{Re}\,\Big( \frac 1{2i\nu} \xi_\bullet e^{i\nu t}
                 \Big) , \qquad
      \xi_\bullet
  := i\nu\xi_0 + \dot\xi_0 , \\
\nonumber
      \xi_1(\nu,t)
  &=& \frac{\kappa}{\rho\nu}
        \int_0^t \sin \nu(t-\tau)\, \psi(\tau)\, d\tau
 \, =\,  \frac{\kappa}{2i\rho\nu}
        \int_0^t
         \big( e^{i\nu t - i\nu\tau} - e^{-i\nu t + i\nu\tau}
         \big)\, \psi(\tau)
         \, d\tau, \\
\label{xi1=}
  &=& \frac{\kappa}{2i\rho\nu}
           \Big( e^{i\nu t} \hat\psi_t(\nu)
                - e^{-i\nu t} \hat\psi_t(-\nu)
           \Big)
   = \mbox{Re}\,
        \Big( \frac{\kappa}{i\rho\nu}
              \hat\psi_t(\nu)\, e^{i\nu t}
        \Big)  .
\end{eqnarray}

Let $\|\cdot\|_{L_2}$ be the $L_2$-norm:
$   \| g \|^2_{L_2}
  = \int |g(\nu)|^2 \, d\nu \quad
    \mbox{for any $g:\mR\to\mC$}$.

\begin{lem}
\label{lem:E=||||}
The function $\EE(t) = \EE(\xi(\nu,t),\dot\xi(\nu,t))$, defined in (\ref{psiE}), satisfies the equation
$$
    \EE(t)
  = \| \hat a\hat\psi_t + \sqrt{\rho}\xi_\bullet \|^2_{L_2}.
$$
\end{lem}

{\it Proof}. By direct computation we obtain:
$\EE =  \EE_{22} + \EE_{21} + \EE_{11}$, where
\begin{eqnarray*}
         \EE_{22}(t)
   &:=&  \int \rho ( \nu^2\xi_2^2 + \dot\xi_2^2 ) \, d\nu
  \equiv  \EE(0)
     =    \| \sqrt{\rho}\,\xi_\bullet \|^2_{L_2}, \\
         \EE_{21}(t)
   &:=&  \int \rho ( \nu^2\xi_2\xi_1 + \dot\xi_2\dot\xi_1 ) \, d\nu
     =   \int \hat a \hat\psi_t \,
          \overline {\sqrt{\rho} \xi_\bullet}
          \, d\nu
      + \int \overline {\hat a \hat\psi_t} \,
          \sqrt{\rho} \xi_\bullet
          \, d\nu     , \\
         \EE_{11}(t)
   &:=&  \int \rho ( \nu^2\xi_1^2 + \dot\xi_1^2 ) \, d\nu
     =   \|\hat a\, \hat\psi_t\|^2_{L_2}.
\end{eqnarray*}
\qed

By Proposition \ref{cor:ene-ineq} $\EE(t)\le c_3(E)$. Therefore we have the following

\begin{cor}
\label{cor:rhoxi}
$\| \sqrt{\rho}\xi_\bullet \|^2_{L_2} \le c_3(E)$ and
$\|\hat a \hat\psi_t \|^2_{L_2} < 2c_3(E)$.
\end{cor}

\section{Spectrum of $\int \kappa\xi\, d\nu$}
\label{sec:kappaxi}

We put
\begin{equation}
\label{ggA}
     \hat\gamma(\nu)
  =  \pi\frac{\kappa\xi_\bullet(\nu)}{i\nu}, \quad
     \AAA_0^+(t)
  =  \int_t^\infty \psi(t - \tau) \, w(\tau)\, d\tau, \quad
     \AAA_0^-(t)
  =  -\int^t_{-\infty} \psi(t - \tau) \, w(\tau)\, d\tau .
\end{equation}

\begin{lem}
\label{lem:hatgamma}
$\hat\gamma\in L_1$.
\end{lem}

\begin{lem}
\label{lem:AAA}
$ \lim_{t\to +\infty} \AAA_0^+(t)
 = \lim_{t\to -\infty} \AAA_0^-(t) = 0$.
\end{lem}

Proofs of Lemmas \protect\ref{lem:hatgamma}--\ref{lem:AAA} are contained in Section \ref{sec:lemmas}.

\begin{prop}
\label{prop:cxi}
Let $\gamma = \FF^{-1}(\hat\gamma)$. Then
\begin{equation}
\label{limlim}
     \int \kappa\xi \, d\nu
  = \FF^{-1}(\hat w_\diamond^\pm\hat\psi)
    + \mbox{\rm Re}
        \big( \gamma(t) + \AAA_0^\pm(t) \big).
\end{equation}
\end{prop}

\begin{cor}
Assertion (2) of Theorem \ref{theo:exotic} holds:
$$
     \int \kappa \xi\, d\nu
  = \FF^{-1}(\hat w_\diamond^\pm \hat\psi)
         + \thet^\pm , \qquad
     \thet^\pm\in \CC_0^\pm .
$$
\end{cor}

{\it Proof of Proposition \ref{prop:cxi}}. First note that
$$
    \int \kappa \xi_2 \, d\nu
  = \frac{1}{2\pi}
     \mbox{Re}\, \int e^{i\lambda t} \hat\gamma(\lambda)\, d\lambda
  =  \mbox{Re}\,\gamma(t).
$$

By (\ref{xi1=})
$$
      \int \kappa\xi_1 \, d\nu
  =   \mbox{Re}\, \AAA(t) , \qquad
      \AAA(t)
 :=   \frac 1{2\pi}
      \int e^{i\nu t} \hat w(\nu) \hat\psi_t(\nu) \, d\nu
  =   \FF^{-1} (\hat w\, \hat\psi_t)(t)
  =   (\psi_t * w)(t) .
$$
By {\bf H5} and Lemma \ref{lem:a-w} $w\in L_1$. Therefore
$\AAA(t) = \AAA_1^\pm(t) - \AAA_0^\pm(t)$, where $\AAA_0^\pm$ are defined by (\ref{ggA}) and
$$
     \AAA_1^\pm(t)
  =  \int \psi(t - \tau) \, \chi_{[0,\pm\infty)}(\tau)
                         \, w(\tau)\, d\tau
  =  \big( \psi * ( \chi_{[0,\pm\infty)} w ) \big)(t).
$$
The function $\psi$ is real while $\hat w$ is purely imaginary. Therefore $\overline w(\tau) = - w(-\tau)$, $\tau\in\mR$. Hence
$$
    \mbox{Re}\, \AAA_1^\pm
  = \frac12
    \psi * \big( \chi_{[0,\pm\infty)}(\tau) w(\tau)
               - \chi_{[0,\pm\infty)}(\tau) w(-\tau) \big)
  = \FF^{-1}\big(\hat w_\diamond^\pm \hat\psi\big),
$$
where we put
\begin{equation}
\label{hatwdiamond=}
     \hat w_\diamond^\pm
  =  \frac12 \FF\Big( \chi_{[0,\pm\infty)} \,
                       \big( w(\tau) - w(-\tau) \big)
                \Big)
  =  \frac 1{4\pi}
        \Big( \hat\chi_{[0,\pm\infty)}
              * \big( \hat w(\lambda) - \hat w(-\lambda) \big)
        \Big).
\end{equation}
By using an explicit form for the distributions
$\hat\chi_{[0,\pm\infty)}$ \cite{Vlad}:
$$
    \hat\chi_{[0,\pm\infty)}(\lambda)
  = \pm\pi\delta(\lambda) + \frac{1}{i}\,\mbox{v.p.}\frac1\lambda ,
$$
we obtain (\ref{wdiamond}).
\qed

\begin{rem}
\label{rem:diamond}
By the first equation (\ref{hatwdiamond=}) if $w$ lies in $L_1$ then $w^\pm_\diamond$ also lie in $L_1$. Since
$\hat w(\lambda) - \hat w(-\lambda)$ is odd and purely imaginary, by (\ref{wdiamond}) $\mbox{\rm Im}\,\hat w^\pm_\diamond$ are odd and
$\mbox{\rm Re}\,\hat w^\pm_\diamond$ are even.
\end{rem}

\section{Three proofs}
\label{sec:lemmas}

{\bf Proof of Lemma \ref{lem:hatgamma}}. Corollary \ref{cor:rhoxi}
and Hypothesis {\bf H2} imply the inclusion
\begin{equation}
\label{inL2}
  \sqrt{\rho} \xi_\bullet,
  \kappa / (\sqrt{\rho} \nu)
  \in L_2.
\end{equation}
Thus $\kappa\xi_\bullet / \nu \in L_1$. \qed

{\bf Proof of Theorem \ref{theo:exotic}, assertion (1)}

{\bf 1}. The family of distributions $\hat a \hat\psi_t \in\hat\CC_b$ tends as $t\to +\infty$ to
$$
     \hat a\hat\psi_+
  := \hat a\,\FF(\chi_{[0,+\infty]} \psi) \in \hat\CC_b .
$$
Indeed, by {\bf H5} $a = \FF^{-1}(\hat a)\in\MM$. Therefore for any
$\hat\ph = \FF(\ph) \in\hat\MM$ we have: $a*\overline\ph\in\MM$. This implies that
$$
    (\hat a\hat\psi_t,\hat\ph)
  = 2\pi \big( \FF^{-1}(\hat a\hat\psi_t), \ph \big)
  = 2\pi \int (a * \psi_t) \, \overline\ph \, d\tau
  = 2\pi \int \psi_t \, (a*\overline\ph)\, d\tau .
$$
By Lemma \ref{lem:N} when $t\to +\infty$ this tends to
$   2\pi \int \psi_+ \, (a*\overline\ph)\, d\tau
  = (\hat a\hat\psi_+,\hat\ph)$.

{\bf 2}. By Corollary \ref{cor:rhoxi} the $L_2$ norm of
$\hat a\hat\psi_+$ as a functional on $\hat\MM\cap L_2 \subset L_2$ does not exceed $\sqrt{2c_3(E)}$. Since $\hat\MM\cap L_2$ is dense in $L_2$, by the Hahn-Banach theorem $\hat a\hat\psi_+$ can be uniquely extended
to a continuous functional on $L_2$ with the same norm. We denote
this functional again by $\hat a\hat\psi_+$. By the Riesz theorem
$\hat a\hat\psi_+$ can be identified with some element of $L_2$.

{\bf 3}. Analogously we have: $\hat a\hat\psi_- \in L_2$, where
$\hat\psi_- = \FF(\chi_{(-\infty,0]} \psi)$. Since
$\hat\psi = \hat\psi_+ + \hat\psi_-$, we have: $\hat a\hat\psi\in L_2$.

By (\ref{inL2})
$\displaystyle{
 \frac{\kappa}{\sqrt{\rho} \nu}
 = \frac{\hat a}{\nu}\in L_2}$. Therefore
$\displaystyle{
 \hat w\hat\psi = \frac{2\pi\hat a}{i\nu} \cdot \hat a\hat\psi
               \in L_1}$.
In particular, $\FF^{-1}(\hat w\hat\psi) \in \CC_0$.

{\bf 4}. Take any $\lambda\in\sigma$. By assertion (a) of Lemma \ref{lem:1/Phi} there exists an interval $I\subset\sigma$,
$\lambda\in I$ such that for any $\hat\mu\in\hat\MM$ with
$\supp\hat\mu\subset I$ we have: $\hat\mu = \hat g\hat w$, where
$\hat g\in\hat\MM$. Then
$$
    \FF^{-1}(\hat\mu\hat\psi)
  = \FF^{-1}(\hat g\hat w\hat\psi)
  = g * \FF^{-1}(\hat w\hat\psi).
$$
Since $g\in\MM$ and $\FF^{-1}(\hat w\hat\psi)\in\CC_0$, by Lemma \ref{lem:C0*L1} $g * \FF^{-1}(\hat w\hat\psi)\in\CC_0$. Hence $\lambda\not\in\ssupp\hat\psi$. \qed

{\bf Proof of Lemma \ref{lem:AAA}}.  We have:
$   |\AAA_0^\pm(t)|
 =  |\psi * (\chi_{[t,\pm\infty)} w)(t)|
\le \|\chi_{[t,\pm\infty)} w\|_\MM\, \|\psi\|_\infty $.
By Corollary \ref{cor:chimu}
$\lim_{t\to\pm\infty} \|\chi_{[t,\pm\infty)} w\|_\MM = 0$. \qed

\section{The spaces $\CC_b$, $\CC_0^\pm$, $\MM$ and their Fourier images}
\label{sec:CCM}

\subsection{Definitions}
\label{subsec:def}

{\bf 1}. (a) Let $\CC_b$ be the space of uniformly continuous bounded functions $\mR\to\mC$. It is a Banach space with respect to the $L_\infty$-norm:
$$
    \|\ph\|_\infty
  = \sup |\ph|, \qquad
    \ph\in\CC_b.
$$
Moreover, $\CC_b$ is an associative commutative Banach algebra: for any $\ph,\psi\in\CC_b$ we have: $\ph\psi\in\CC_b$ and
$\|\ph\psi\|_\infty\le \|\ph\|_\infty \,\|\psi\|_\infty$.

(b) Let $\CC_0^+\subset\CC_b$ be the space of functions vanishing at $+\infty$ and $\CC_0^-\subset\CC_b$ the space of functions vanishing at $-\infty$. Then $\CC_0^\pm$ are closed vector subspaces and moreover, ideals in $(\CC_b,\cdot)$, i.e., $\CC_0^\pm\cdot\CC_b\subset\CC_0^\pm$.

We define $\CC_0 = \CC_0^+\cap\CC_0^-$. It is also a closed ideal in $\CC_b$.
\medskip

(c) Let $\MM$ be the space of densities of finite Radon $\mC$-measures on $\mR$. These densities are distributions on $\mR$, so for any $\mu\in\MM$ the corresponding measure is $\mu(\tau)\, d\tau$. The space $\MM$ is conjugated to $\CC_b$. Therefore it is a Banach space with respect to the norm
$$
    \|\mu\|_{\MM}
  = \sup_{\|\ph\|_\infty = 1} \mu(\ph), \qquad
    \mu(\ph) = \int \ph\mu\, d\tau, \quad
    \ph\in\CC_b, \; \mu\in\MM.
$$
Obviously, $L_1\subset\MM$.

Recall that $\MM$ is a convolution Banach algebra:
$\|\mu_1 * \mu_2\|_{\MM} \le \|\mu_1\|_{\MM}\,\|\mu_2\|_{\MM}$.
Any $\mu\in\MM$ admits an expansion
\begin{equation}
\label{expansion}
   \mu = \mu_1 - \mu_2 + i\mu_3 - i\mu_4, \quad
   \mbox{where $\mu_1,\mu_2,\mu_3,\mu_4\in\MM$ are non-negative},
\end{equation}
see for example, \cite{Boga}.

(d) We also define the convolution of $\mu\in\MM$ and $\ph\in\CC_b$.
Consider the linear isometric operator
$$
  \bs_t : \CC_b\to\CC_b, \quad
  (\bs_t\ph)(\tau) = \ph(t - \tau).
$$
Then $(\mu * \ph)(t) = \mu (\bs_t\ph)$.
By Lemma \ref{lem:C0*L1} we have:
$\MM *\CC_b\subset\CC_b$, $\MM *\CC_0\subset\CC_0$ and moreover,
$$
      \|\mu * \ph\|_\infty
  \le \|\mu\|_{\MM} \, \|\ph\|_\infty.
$$
\medskip

\noindent
{\bf 2}. (a) We put $\hat\CC_b = \FF(\CC_b)$. Elements of the space $\hat\CC_b$ are distributions. In particular, $\hat\CC_b$ contains $\delta$-function. The equation
$\FF^{-1}(\hat\ph *\hat\psi) = 2\pi\ph\psi$ implies that $(\hat\CC_b,*)$ is an associative commutative algebra: for any $\hat\ph,\hat\psi\in\hat\CC_b$ we have:
$\hat\ph * \hat\psi\in\hat\CC_b$.

(b) Moreover, $\hat\CC_0^+ = \FF(\CC_0^+)$,
$\hat\CC_0^- = \FF(\CC_0^-)$, and $\hat\CC_0 = \FF(\CC_0)$ are ideals in $(\hat\CC_b,*)$.

(c) We also put $\hat\MM = \FF(\MM)$. It is easy to show that $\hat\MM\subset\CC_b$, however $\hat\MM\ne\CC_b$, \cite{Boga}. By the Riemann-Lebesgue theorem the Fourier image of $L_1\subset\MM$ lies in $\CC_0$.

By {\bf 1}(c) for any $\hat\mu_1,\hat\mu_2\in\hat\MM$ we have: $\hat\mu_1\hat\mu_2\in\hat\MM$.

(d) By {\bf 1}(d) we have: $\hat\MM\cdot\hat\CC_b\subset\hat\CC_b$,
$\hat\MM\cdot\hat\CC_0^+\subset\hat\CC_0^+$,
$\hat\MM\cdot\hat\CC_0^-\subset\hat\CC_0^-$, and
$\hat\MM\cdot\hat\CC_0\subset\hat\CC_0$.
\medskip

\noindent
{\bf 3}. (a) Let $\hat\ph\in\hat\CC_b$, $\lambda\in\mR$. We say that $\lambda\not\in\ssupp^+(\hat\ph)$ iff there exists an interval
$I\ni\lambda$ such that for any $\hat\mu\in\hat\MM$ with
$\supp\hat\mu\subset I$ we have: $\hat\mu\hat\ph\in \hat\CC_0^+$.

(b) Equivalently, for any $\hat\ph\in\hat\CC_b$ we say that $\lambda\in\ssupp^+(\hat\ph)$ iff for any interval
$I\ni\lambda$ there is a function $\hat\mu\in\hat\MM$ such that
$\supp\hat\mu\subset I$ and $\hat\mu\hat\ph\in\hat\CC_b\setminus\hat\CC_0^+$.

Analogously we define $\ssupp^-(\ph)$ and $\ssupp(\ph)$.
The sets $\ssupp^+(\hat\ph)$, $\ssupp^-(\hat\ph)$, and
$\ssupp(\hat\ph) = \ssupp^+(\hat\ph) \cup \ssupp^-(\hat\ph)$ are obviously closed.

\subsection{Lemmas on $\CC_b$, $\CC_0^\pm$, $\MM$ and their Fourier images}

We start from some notation. Consider the functions
$$
     h_0(\tau)
  =  \left\{ \begin{array}{cl}
              1 - |\tau|\quad & \mbox{if } |\tau| \le 1, \\
                  0     \quad & \mbox{if } |\tau| \ge 1.
             \end{array}
     \right. \quad
     h_k(\tau)
  =  h_0(k - \tau), \qquad
     k\in\mZ .
$$
Obviously, $h_k\in\CC_0$ form a partition of unity:
$\sum_{k\in\mZ} h_k = 1$.

For any $K\in\mN$ we put
$$
     \bh_K^<
  =  \sum_{|k|\le K} h_k, \quad
     \bh_K^>
  =  \sum_{|k| > K} h_k.
$$

\begin{lem}
\label{lem:N}
Let $\mu\in\MM$. Then for any $\eps > 0$ there exists $N>0$ such that for any $\ph\in\CC_b$ such that $\supp\ph\cap [-N,N] = \emptyset$ we have: $|\mu(\ph)| < \eps\,\|\ph\|_\infty$.
\end{lem}

{\it Proof}. By using the expansion (\ref{expansion}), it is sufficient to assume that $\mu$ is non-negative. Then
$$
    \|\mu\|_{\MM}
  = \mu(1)
  = \sum_{k\in\mZ} \mu(h_k).
$$
This series converges. Therefore for any $\eps>0$ there exists $K\in\mN$ such that $\mu(\bh_K^>) < \eps$. For any $\ph\in\CC_b$ such that
$\supp\ph \cap [-K-2,K+2] = \emptyset$ we have:
$$
     |\mu(\ph)|
  =  |\mu(\bh_K^>\,\ph)|
 \le \|\ph\|_\infty \, |\mu(\bh_K^>)|
  <  \eps\, \|\ph\|_\infty .
$$
\qed

\begin{cor}
\label{cor:chimu}
For any $\mu\in\MM$ we have:
$\lim_{t\to +\infty} \|\chi_{[t,+\infty)}\,\mu\|_{\MM}
 = \lim_{t\to -\infty} \|\chi_{(-\infty,t]}\,\mu\|_{\MM} = 0$.
\end{cor}

Indeed, since $\mu$ can be assumed to be non-negative,
$$
     \|\chi_{[t,+\infty)}\,\mu\|_{\MM}
  =  (\chi_{[t,+\infty)}\,\mu)(1)
 \le \mu(\tilde h), \qquad
     \tilde h
  =  \sum_{k = [t]-1}^{+\infty} h_k.
$$
Now it is sufficient to note that
$\supp(\tilde h) \cap (-t+2,t-2) = \emptyset$,
$\|\tilde h\|_\infty = 1$. Then by Lemma \ref{lem:N}
$\lim_{t\to +\infty} \|\chi_{[t,+\infty)}\,\mu\|_{\MM} = 0$.

The second limit can be considered analogously. \qed

\begin{lem}
\label{lem:C0*L1}
For any $\ph\in\CC_b$ and $\mu\in\MM$ we have: $\mu *\ph\in\CC_b$ and
$\|\mu*\ph\|_\infty \le \|\mu\|_{\MM} \, \|\ph\|_\infty$.

Moreover, if $\ph\in\CC_0^\star$, $\star\in\{+,-\}$, then $\mu*\ph\in\CC_0^\star$.
\end{lem}

{\it Proof}. The function $\mu*\ph$ is bounded because
$$
      |(\mu*\ph)(t)|
   =  |\mu (\bs_t\ph|
  \le  \|\mu\|_{\MM}\,\|\ph\|_\infty.
$$
Since $\ph$ is uniformly continuous, for any $\eps > 0$ there exists an $\alpha > 0$ such that for any $t$ we have:
$|\bs_{t+\alpha}\ph - \bs_t\ph| < \eps$. Then
$$
       |(\mu*\ph)(t+\alpha) - (\mu*\ph)(t)|
   =  |\mu (\bs_{t+\alpha}\ph - \bs_t\ph|
  \le \|\mu\|_{\MM}\, \eps.
$$
This implies uniform continuity of $\ph * h$.

Now suppose that $\ph\in\CC_0^+$. Take an arbitrary $\eps>0$. By Lemma \ref{lem:N} there exists $N > 1$ such that
\begin{equation}
\label{mupsi<eps}
  |\mu(\psi)| < \eps \|\psi\|_\infty \quad
  \mbox{for any $\psi\in\CC_b$ such that }
  \supp\psi\cap (-N-1,N+1) = \emptyset.
\end{equation}
If necessary we take a larger $N$ such that
\begin{equation}
\label{ph<eps}
    |\ph(\tau)| < \eps \quad
    \mbox{for any $\tau > N$}.
\end{equation}
We have:
$(\mu*\ph)(t) = \mu(\bs_t\ph)
              = \mu(\bh_N^<\, \bs_t\ph) + \mu(\bh_N^>\, \bs_t\ph)$.
We estimate the last two terms separately.

For $t > 2N$ by (\ref{ph<eps}) $\|\bh_N^<\, \bs_t\ph\|_\infty < \eps$. Therefore $|\mu(\bh_N^<\, \bs_t\ph)| < \eps \|\mu\|_{\MM}$.

On the other hand $\supp(\bh_N^>\,\bs_t\ph)\cap (-N-1,N+1) = \emptyset$. Hence,
$|\mu(\bh_N^>\, \bs_t\ph)| < \eps \|\bh_N^>\, \bs_t\ph\|_\infty
                        \le \eps \|\ph\|_\infty$.

Combining these two estimates, we have:
$|(\mu*\ph)(t)| \le \eps (\|\mu\|_{\MM} + \|\ph\|_\infty)$.

The case $\ph\in\CC_0^-$ is analogous.
\qed

\begin{cor}
\label{cor:C0*L1}
For any $\hat\ph\in\hat\CC_b$, $\hat\ph_0\in\hat\CC_0^\star$, $\star\in\{+,-\}$, and $\hat\mu\in\hat\MM$ we have: $\hat\mu\hat\ph\in\hat\CC_b$ and
$\hat\mu\hat\ph_0\in\hat\CC_0^\star$.
\end{cor}

\begin{lem}
\label{lem:delta}
Let $\hat\beta\in\hat\MM$, $\hat\beta(0) = 1$, and let $\hat\beta_r(\lambda)=\hat\beta(r\lambda)$. Then for any $\hat\ph\in\hat\CC_b$
$$
  \lim_{r\to 0} \big\| \FF^{-1}(\hat\beta_r\hat\ph) - \ph
                \big\|_\infty
  = 0.
$$
\end{lem}

{\it Proof}. First note that
$$
  \beta_r(\ph) = \beta(\ph_r) \quad
  \mbox{for any $\ph\in\CC_b$, where }
  \ph_r(\tau) = \ph(r\tau).
$$
Take any $\eps > 0$. There exists $b>0$ such that
\begin{equation}
\label{ph-ph}
  |\ph(t-\tau) - \ph(t)| < \eps \quad
  \mbox{for any $t\in\mR$ and $|\tau| \le b$}.
\end{equation}
By Lemma \ref{lem:N} we can choose $r_0 > 0$ such that
\begin{equation}
\label{betaph}
    |\beta(\ph)|
  < \eps \|\ph\|_\infty \quad
    \mbox{for any $\ph\in\CC_b$ with }
    \supp(\ph) \cap [2-b/r_0,-2+b/r_0] = \emptyset.
\end{equation}
Below we assume that $0 < r < r_0$.
Choosing an integer $K\in [-2+b/r,-1+b/r]$, we obtain
$$
    (\beta_r * \ph)(t)
  = \beta_r(\bs_t\ph)
  = \beta(\bs_{t/r}\ph_r)
  = \ph(t)\,\beta(\bh_K^<)
   + \beta\big(\bh_K^< (\bs_{t/r}\ph_r - \ph(t)) \big)
   + \beta(\bh_K^> \, \bs_{t/r}\ph_r).
$$
We have the following estimates.

(1) Since $\beta(1) = 1$, and
\begin{equation}
\label{supp(h)}
  \supp(\bh_K^>)\cap [-2-b/r,2+b/r] = \emptyset,
\end{equation}
the inequality (\ref{betaph}) implies that
$$
     |\beta(\bh_K^<) - 1|
  =  |\beta(\bh_K^>)|
 \le \eps\, \|\bh_K^>\|_\infty = \eps.
$$

(2) Since by (\ref{ph-ph})
$$
  |(\bs_{t/r}\ph_r)(\tau) - \ph(t)| < \eps \quad
  \mbox{for any $|\tau| \le b/r$},
$$
the inequality $K \le b/r_0 - 1$ implies
$\big\|\bh_K^< (\bs_{t/r}\ph_r - \ph(t)) \big\|_\infty < \eps$,
we have:
$$
      \big| \beta\big(\bh_K^< (\bs_{t/r}\ph_r - \ph(t)) \big) \big|
  \le \|\beta\|_{\MM}\,
      \big\|\bh_K^< (\bs_{t/r}\ph_r - \ph(t)) \big\|_\infty
  \le \eps \, \|\beta\|_{\MM} .
$$

(3) Again by using (\ref{supp(h)}), we obtain:
$$
      |\beta(\bh_K^> \, s_{t/r}\ph_r)|
  \le \eps \|\ph\|_\infty .
$$

Combining estimates (1)--(3), we obtain:
$|(\beta_r * \ph)(t) - \ph(t)|
 < \eps\, ( 2\|\ph\|_\infty + \|\beta\|_{\MM})$.
\qed

\begin{lem}
\label{lem:1/Phi}
Suppose that $\hat\mu\in\hat\MM$, $\hat\mu(\lambda)\ne 0$, and $I$ is a sufficiently small interval containing $\lambda$. Then

(a) for any $\hat\mu_0\in\hat\MM$, $\supp\hat\mu_0\subset I$
  we have $\hat\mu_0 = \hat\mu \hat\mu_1$ for some $\hat\mu_1\in\hat\MM$,

(b) for any $\hat\ph\in\hat\CC_b$, $\supp\hat\ph\subset I$
  we have $\hat\ph = \hat\mu\hat\ph_1$ for some $\hat\ph_1\in\hat\CC_b$.
\end{lem}

{\it Proof}. Without loss of generality we have: $\lambda=0$. We assume that $\supp\hat\mu_0 = [-\alpha,\alpha]\subset I$. Below $\alpha$ will be chosen sufficiently small. We put
$\beta = \chi_{[-1,1]}$ and $\beta_\eps(\tau) = \beta(\eps\tau)$. As usual, $\hat\beta_\eps = \FF(\beta_\eps)$. Then
$$
    \frac{\hat\mu_0}{\hat\mu}
  = \frac{\hat\mu_0}
         {\hat\beta_\eps *\hat\mu
           - (\hat\mu - \hat\beta_\eps *\hat\mu)}
  = \frac{\hat\mu_0}{\hat\beta_\eps *\hat\mu}
   + \frac{\hat\mu_0}{(\hat\beta_\eps *\hat\mu)^2}
       \big( \hat\mu - \hat\beta_\eps *\hat\mu \big)
   + \frac{\hat\mu_0}{(\hat\beta_\eps *\hat\mu)^3}
       \big( \hat\mu - \hat\beta_\eps *\hat\mu \big)^2
           + \ldots
$$
Let $\hat h$ be a smooth function such that
$\hat h|_{[-1,1]} = 1$ and $\supp(\hat h) \subset [-2,2]$. Then
putting
$$
    \hat h_\alpha(\lambda)
  = \hat h(\lambda / \alpha), \quad
    \hat f_{\eps,\alpha}
  = \hat h_\alpha / (\hat\beta_\eps * \hat\mu), \quad
    f_{\eps,\alpha}
  = \FF^{-1}(\hat f_{\eps,\alpha}), \quad
    g_\eps
  = \FF^{-1}(\hat\mu - \hat\beta_\eps * \hat\mu),
$$
we have: $\hat h|_{\supp(\hat\mu_0)} = 1$. Therefore
$$
      \FF^{-1}\Big( \frac{\hat\mu_0}{\hat\mu} \Big)
   =  \sum_{k=0}^\infty
       \FF^{-1}\Big( \frac{\hat\mu_0 (\hat h_\alpha)^{k+1}}
                          {(\hat\beta_\eps *\hat\mu)^{k+1}}
               \Big)
      * \FF^{-1}\Big(
                   \big( \hat\mu - \hat\beta_\eps *\hat\mu
                   \big)^k \Big)
   =  \sum_{k=0}^\infty
          \mu_0 * f_{\eps,\alpha}^{*k}
      * g_\eps^{*k} .
$$

If the interval $I$ and $\eps > 0$ are sufficiently small,
$(\hat\beta_\eps * \hat\mu)|_I$ is close to $\hat\mu(0)\ne 0$.
Moreover, since $\hat\beta_\eps$ and $\hat h_\alpha$ are smooth,
$\hat f_{\eps,\alpha}$ is also smooth, and $f_{\eps,\alpha}$ is fast decreasing. Since $\supp(\hat h_\alpha)$ is compact, $f_{\eps,\alpha}$ is smooth. The following equations hold
\begin{equation}
\label{two_limits}
    \lim_{\eps\to+0} \|g_\eps\|_{\MM}
  = 0 \qquad
    \mbox{and for any $\eps > 0$} \quad
    \lim_{\alpha\to+0} \|f_{\eps,\alpha}\|_{\MM}
  = |b_\eps|\, \|h\|_{\MM} ,
\end{equation}
where we put $1 / b_\eps = (\hat\beta_\eps * \hat\mu)(0)$.
We will prove them a few lines below.

By (\ref{two_limits}) we can choose $\eps,\alpha > 0$ so that
$\|f_{\eps,\alpha}\|_{\MM} \, \|g_\eps\|_{\MM} \le 1/2$. The estimate
$$
      \|\FF^{-1}(\hat\mu_0 / \hat\mu)\|_{\MM}
  \le \sum_{k=0}^\infty \|\mu_0\|_{\MM}\,
                     \|f_{\eps,\alpha}\|_{\MM}^k \,
                  \|g_\eps\|_{\MM}^k
   <  \infty
$$
completes the proof of assertion (a).

Proof of assertion (b) is analogous: the final estimate is
$$
      \|\FF^{-1}(\hat\ph / \hat\mu)\|_\infty
  \le \sum_{k=0}^\infty \|\ph\|_\infty\,
                     \|f_{\eps,\alpha}\|_{\MM}^k \,
                  \|g_\eps\|_{\MM}^k
   <  \infty .
$$

Finally we prove equations (\ref{two_limits}). The first one follows from Corollary \ref{cor:chimu}. To prove the second one we note that
$\hat h_\alpha = \hat h_\alpha \hat h_{2\alpha}$ for any $\alpha > 0$.
Then
$$
    \hat f_{\eps,\alpha}
  = b_\eps \hat h_\alpha
   + J \hat h_\alpha, \qquad
    J
  = \hat h_{2\alpha}
       \Big( \frac{1}{\hat\beta_\eps * \hat\mu} - b_\eps \Big).
$$
Hence
$    \|\hat f_{\eps,\alpha}\|_{\MM}
 \le \|b_\eps \hat h_\alpha\|_{\MM}
   + \|J\|_{\MM} \|\hat h_\alpha\|_{\MM}.
$
Now it remains to note that $\|h_\alpha\|_{\MM} = \|h\|_{\MM}$ for all positive $\alpha$ and $\lim_{\alpha\to+0} \|J\|_{\MM} = 0$ because the function $1 / (\hat\beta_\eps * \hat\mu) - b_\eps$ vanishes at zero, is independent of $\alpha$ and smooth in a small neighborhood of zero.
\qed

\begin{lem}
\label{lem:a-w}
For any $\hat a\in\hat\MM$ such that $\hat a'\in\hat\MM$ and
$\hat a(0) = 0$ the Fourier transform of $\hat a^2(\lambda) / \lambda$ lies in $L_1$.
\end{lem}

{\it Proof}. We put $\ba = \FF^{-1}(\hat a / \lambda)$. Then
$\ba = i\,\theta * a$, where
$\theta = \frac12(\chi_{(0,+\infty)} - \chi_{(-\infty,0)})$.
Since $a\in\MM$, $\theta * a$ is defined almost everywhere and bounded. Let us show that its $L_1$-norm (or $\MM$-norm, what is the same) is finite. Indeed, since $\hat a(0) = 0$, for almost all $s\in\mR$ we have:
$$
    \int_{-\infty}^s a(\tau)\, d\tau
  = - \int_s^{+\infty} a(\tau)\, d\tau .
$$
Then
\begin{eqnarray*}
       \|\theta * a\|_{\MM}
  &=&   \int ds \,\Big| \int_{-\infty}^s a(\tau)\, d\tau
                       - \int_s^{+\infty} a(\tau)\, d\tau
                 \Big| \\
  &=&   2\int_{-\infty}^0 ds \,
             \Big| \int_{-\infty}^s a(\tau)\, d\tau \Big|
      + 2\int_0^{+\infty} ds \,
             \Big| \int_s^{+\infty} a(\tau)\, d\tau \Big| \\
 &\le&  2\int_{-\infty}^0 ds \,
              \int_{-\infty}^s |a(\tau)|\, d\tau
      + 2\int_0^{+\infty} ds \,
              \int_s^{+\infty} |a(\tau)|\, d\tau  \\
  &=&   2\int |\tau\, a(\tau)| \, ds
   =    2 \|\FF^{-1}(\hat a')\|_{\MM}
   <   \infty .
\end{eqnarray*}

We have: $\FF^{-1}(\hat a^2 / \lambda) = \ba * a$. Since the function $\ba$ is bounded and $a\in\MM$ then $\ba * a$ is also a bounded function. The estimate
$$
      \|\ba * a\|_{\MM}
  \le \|\ba\|_{\MM}\, \|a\|_{\MM} < \infty .
$$
completes the proof. \qed

\subsection{Lemmas on the singular support}

\begin{lem}
\label{lem:loc->glob}
Suppose that $\hat\ph\in\hat\CC_b$, $-\infty\le a < b\le +\infty$, and
$(a,b)\subset \mR\setminus \ssupp^\star(\hat\ph)$, $\star\in\{+,-\}$. Then for any $\hat\mu\in\hat\MM$, where
$\mbox{\rm supp}(\hat\mu)\subset [a,b]$, we have:
$\hat\mu\hat\ph \in\hat\CC_0^\star$.
\end{lem}

\begin{figure}
\begin{center}
\includegraphics[scale=1.3]{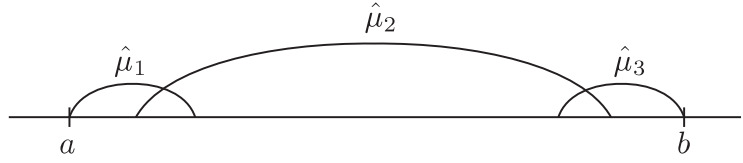}
\caption{Graphs of the functions $\hat\mu_1,\hat\mu_2$,
         and $\hat\mu_3$.}
\label{fig:hhh}
\end{center}
\end{figure}

{\it Proof}. For definiteness we take $\star = +$. First suppose that $a$ and $b$ are finite. Take an arbitrary $\eps > 0$. We can represent $\hat\mu$ in the form $\hat\mu = \hat\mu_1 + \hat\mu_2 + \hat\mu_3$, where $\hat\mu_j\in\hat\MM$,
$$
  \mbox{supp}(\hat\mu_1) \subset [a,b],\quad
  \mbox{supp}(\hat\mu_2) \subset (a,b),\quad
  \mbox{supp}(\hat\mu_3) \subset [a,b],\quad
  \mbox{and}\quad
    \|\mu_1\|_{\MM},\|\mu_3\|_{\MM}
  < \frac12 \eps / \|\ph\|_\infty,
$$
see Fig. \ref{fig:hhh}.

Then $\FF^{-1}(\hat\mu\hat\ph) = (\mu_1 + \mu_2 + \mu_3) * \ph$. We have:
$|(\mu_1 + \mu_3) * \ph|
 \le (\|\mu_1\|_{\MM} + \|\mu_3\|_{\MM}) |\ph|_\infty < \eps$
while $\mu_2 * \ph \in\CC_0^+$ because
$\supp\hat\mu_2$ is compact and does not intersect with $\ssupp^+(\hat\ph)$.

Now suppose that $a = -\infty$ and/or $b = +\infty$. Let
$$
  \hat\beta\in\hat\MM, \quad
  \hat\beta(0) = 1,\quad
  \supp(\hat\beta) \subset [-1,1],
$$
and let $\hat\beta_r(\lambda)=\hat\beta(r\lambda)$.
Consider an arbitrary $\eps > 0$ and $\hat\mu\in\hat\MM$ with
$\mbox{\rm supp}(\hat\mu)\subset [a,b]$. Then by Lemma \ref{lem:delta}
for sufficiently small $r$
$$
  \big\| \FF^{-1}\big( \hat\beta_r\hat\mu\hat\ph \big)
        - \FF^{-1}\big(\hat\mu\hat\ph \big)
  \big\|_\infty < \eps.
$$
The set $\supp(\hat\beta_r\hat\mu\hat\ph) \subset [-1/r,1/r]\cap [a,b]$, is compact. Hence $\hat\beta_r\hat\mu\hat\ph\in\hat\CC_0^+$.
Therefore for sufficiently large $\tau > 0$
$$
      |\FF^{-1}(\hat\mu\hat\ph)(\tau)|
  \le \big\| \FF^{-1}\big( \hat\beta_r\hat\mu\hat\ph \big)
        - \FF^{-1}\big(\hat\mu\hat\ph \big)
      \big\|_\infty
     +  \big| (\hat\beta_r\hat\mu\hat\ph)(\tau) \big|
  \le \eps.
$$
The case $\star = -$ is analogous.
\qed

\begin{cor}
\label{cor:ssupp=0}
Suppose that $\hat\ph\in\hat\CC_b$, $\star\in\{+,-\}$, and
$\ssupp^\star(\hat\ph) = \emptyset$. Then $\hat\ph\in\hat\CC_0^\star$. Indeed, it is sufficient to take in Lemma \ref{lem:loc->glob}
$\hat\mu = 1$.
\end{cor}

\begin{lem}
\label{cor:onepoint}
Let $\hat\ph\in\hat\CC_b$, $\star\in\{+,-\}$, and $\ssupp^\star(\hat\ph)$ consist of a finite number of points
$\lambda_1 < \ldots < \lambda_k$. Then
$\hat\ph = \hat\ph_0 + \sum_{j=1}^k \hat\ph_j$, where $\supp\hat\ph_j$ lie in small neighborhoods of $\lambda_j$,
$$
  \hat\ph_0\in\hat\CC_0^\star, \quad
  \hat\ph_j\in\hat\CC_b, \quad
  \ph_j\in C^\infty, \quad
  \Big(\frac{\partial}{\partial\tau} - i\lambda_j\Big) \ph_j
  \in \CC_0^\star.
$$
\end{lem}

{\it Proof}. It is sufficient to take $\hat\ph_j = \hat\mu_j\hat\ph$, where $\hat\mu_j\in\hat\MM$ and for sufficiently small $\eps > 0$
$$
  \supp\hat\mu_j \subset (\lambda_j - 2\eps,\lambda_j + 2\eps), \quad
  \hat\mu_j |_{(\lambda_j - \eps,\lambda_j + \eps)} = 1.
$$
Then $\supp\hat\ph_j$ is compact. Therefore $\ph_j\in C^\infty$.

For any $j = 1,\ldots,n$ $\ssupp^\star(\hat\ph_j) = \{\lambda_j\}$. Let us put
$$
    g_j^+(\lambda)
  = i(\lambda - \lambda_j)\, \chi_{[\lambda_j,+\infty)}(\lambda) ,
    \quad
    g_j^-(\lambda)
  = i(\lambda - \lambda_j)\, \chi_{(-\infty,\lambda_j]}(\lambda) .
$$
Obviously $g^\pm_j\hat\mu_j \in\hat\MM$. Then
$$
    i(\lambda - \lambda_j) \hat\ph_j
  = g_j^+ \hat\mu_j \hat\ph + g_j^- \hat\mu_j \hat\ph,
$$
where by Lemma \ref{lem:loc->glob} $g^\pm_j\hat\mu_j\hat\ph\in\hat\CC_0^\star$. This implies that
$i(\lambda - \lambda_j)\hat\ph_j\in\hat\CC_0^\star$. Hence
$$
    \FF^{-1}(i(\lambda - \lambda_j)\hat\ph_j)
  = \big( \partial / \partial\tau - i\lambda_j \big) \ph_j
 \in \CC_0^\star.
$$
Finally by Corollary \ref{cor:ssupp=0}
$\hat\ph_0 = \hat\ph - \sum_1^k \hat\ph_j \in\hat\CC_0^\star$. \qed

\begin{lem}
\label{lem:ph*psi}
Suppose that $\ssupp^\star(\hat\ph) \subset (-\infty,a]$ and
$\ssupp^\star(\hat\psi) \subset (-\infty,b]$, $\star\in\{+,-\}$.

Then $\ssupp^\star(\hat\ph * \hat\psi) \subset (-\infty,a+b]$.
\end{lem}

{\it Proof}. For any $\eps > 0$ we have:
$\hat\ph = \hat\ph_- + \hat\ph_+$ and
$\hat\psi = \hat\psi_- + \hat\psi_+$, where
\begin{eqnarray*}
 & \hat\ph_-,\hat\psi_- \in \hat\CC_b, \quad
   \hat\ph_+,\hat\psi_+ \in \hat\CC_0^\star, & \\
 & \supp\hat\ph_- \subset (-\infty,a+2\eps], \;
   \supp\hat\psi_- \subset (-\infty,b+2\eps], \;
   \supp\hat\ph_+ \subset [a+\eps,+\infty), \;
   \supp\hat\psi_+ \subset [b+\eps,+\infty).
\end{eqnarray*}
Indeed, it is sufficient to take
$$
    h_+
  = \sum_{k=1}^\infty h_k, \quad
    \hat\ph_+
  = h_+\big( (\tau-a-\eps) / \eps \big)\,\hat\ph, \quad
    \hat\psi_+
  = h_+\big( (\tau-b-\eps) / \eps \big)\,\hat\psi.
$$

Then
$\hat\ph * \hat\psi = \hat\ph_- * \hat\psi_- \bmod\hat\CC_0^\star$,
and $\supp (\hat\ph_- * \hat\psi_-) \subset (-\infty, a+b+4\eps]$.
\qed

\begin{lem}
\label{cor:1/Phi}
Suppose that $\hat\psi\in\hat\CC_b$ satisfies the equation $\hat\Phi\hat\psi = \hat\thet$, where $\hat\thet\in\hat\CC_0^\star$, $\star\in\{+,-\}$, and $\hat\Phi\in\hat\MM$.
Then $\ssupp^\star (\hat\psi) \subset \{\hat\Phi = 0\}$.
\end{lem}

Indeed, let $\hat\Phi(\lambda)\ne 0$. Take a small interval
$I\ni\lambda = 0$ and $\hat\mu_0\in\hat\MM$ with
$\supp\hat\mu_0\subset I$. Then by assertion (b) of Lemma \ref{lem:1/Phi} $\hat\mu_0 / \hat\Phi \in \hat\MM$. Then by Lemma \ref{lem:C0*L1}
$\hat\mu_0\hat\psi = \hat\mu_0 \hat\thet / \hat\Phi \in\hat\CC_0^\star$. \qed

{\bf Acknowledgements}. I am grateful to I.Volovich who attracted my attention to the problem of interaction of an oscillator and thermostat. When working on the paper I had many useful discussions with S.Bolotin, Yu.Drozhinov, A.Shkalikov, and B.Zavialov. The work is partially supported by the grant RFBR 08-01-00681-a and the RAS program ``Mathematical theory of control''.


\begin{thebibliography}{00}
\bibitem{Boga} Bogachev V.I. Measure theory. Vol. I, Springer-Verlag, Berlin, 2007.
\bibitem{Bog} Bogolyubov N.N. An elementary example of a statistical equilibration in a system connected to a thermostat.
    In the book: On some statistical methods in mathematical physics.
    AN USSR publishing 1945, p. 115--137 (in Russian).
\bibitem{Kom} Komech A. Stabilization of the interaction of a string with a nonlinear oscillator. Mosc. Univ. Math. Bull. 46 (6) (1991), 34--39.
\bibitem{KSK97} Komech A. Spohn H. Kunze M.
    Long-time asymptotics for a classical particle interacting with a scalar wave field. Comm. Partial Diff. Eq. 22 (1--2) (1997) 307--335.
\bibitem{KS00} Komech A. Spohn H.
    Long-time asymptotics for a coupled Maxwell-Lorentz equations.
    Comm. Partial Diff. Eq. 25 (3--4) (2000) 559--584.
\bibitem{Riesz-SN}  Riesz F., Sz.-Nagy B.  Functional analysis.
    Dover Books on Advanced Mathematics. Dover Publications, Inc., New York, 1990.
\bibitem{Schwa} Schwartz L. Cours d'analyse. 2. Hermann, Paris, 1981.
\bibitem{Vlad} Vladimirov V.S. Generalized functions in mathematical physics. ``Mir'', Moscow, 1979.
\end{thebibliography}
\end{document}